\documentclass[xcolor=dvipsnames,svgnames,table,reqno]{amsart}

\input xy
\xyoption{all}
\usepackage{accents}
\usepackage{epsfig}
\usepackage{xcolor}
\usepackage{amsthm}
\usepackage{amssymb}
\usepackage{bbm}
\usepackage{amsmath}
\usepackage{amscd}
\usepackage{amsopn}
\usepackage{graphicx}
\usepackage{xspace}

\usepackage{hhline}
\usepackage{easybmat}
\usepackage{caption}   
\usepackage{relsize}

\usepackage{url}
\usepackage{enumitem, hyperref}\hypersetup{colorlinks}

\usepackage{stmaryrd}

\usepackage[T1]{fontenc}

\colorlet{purpleB70}{blue!70!red}

\colorlet{orangeR65}{red!65!yellow}

\definecolor{red2}{HTML}{d41173}

\definecolor{neongreen}{HTML}{1bf702}

\definecolor{radicalred}{HTML}{FF355E}

\definecolor{denim}{HTML}{1560BD}

\definecolor{darkcyan}{rgb}{0.0, 0.55, 0.55}

\definecolor{cilek}{HTML}{FF43A4}

\definecolor{mor}{HTML}{9F00C5}

\definecolor{phlox}{rgb}{0.87, 0.0, 1.0}

\definecolor{fluorescentpink}{HTML}{FF1493}

\definecolor{napiergreen}{rgb}{0.16, 0.5, 0.0}

\definecolor{kellygreen}{rgb}{0.3, 0.73, 0.09}

\definecolor{parisgreen}{HTML}{ 50C878 }

\definecolor{palatinateblue}{rgb}{0.15, 0.23, 0.89}

\definecolor{ceruleanblue}{rgb}{0.16, 0.32, 0.75}

\definecolor{brandeisblue}{rgb}{0.0, 0.44, 1.0}

\definecolor{KLMblue}{HTML}{0FC0FC}

\definecolor{cinnamon}{rgb}{0.82, 0.41, 0.12}

\definecolor{darkorange}{rgb}{1.0, 0.55, 0.0}

\definecolor{darktangerine}{rgb}{1.0, 0.66, 0.07}

\definecolor{deepcarrotorange}{rgb}{0.91, 0.41, 0.17}

\definecolor{internationalorange}{HTML}{FF4F00}

\definecolor{persimmon}{HTML}{EC5800}

\definecolor{pumpkin}{HTML}{FF7518}

\definecolor{darkred}{rgb}{1,0,0} 
\definecolor{darkgreen}{rgb}{0,0.7,0}
\definecolor{darkblue}{rgb}{0,0,1}

\hypersetup{colorlinks,
linkcolor=darkblue,
filecolor=darkgreen,
urlcolor=darkred,
citecolor=darkgreen}

\makeatletter
\def\reflb#1#2{\begingroup
    #2%
    \def\@currentlabel{#2}%
    \phantomsection\label{#1}\endgroup
}
\makeatother

\numberwithin{equation}{section}
\newtheorem{Theorem}{Theorem}
\numberwithin{Theorem}{section}

\newtheorem   {Corollary}[Theorem]{Corollary}
\theoremstyle {definition}

\theoremstyle {remark}
\newtheorem   {Remark}[Theorem]{Remark}
\newtheorem   {Example}[Theorem]{Example}

\newcommand{\FF}{{\mathcal F}}

\newcommand{\Ham}{{\mathit{Ham}}}
\newcommand{\Spec}{{\mathit{Spec}}}

\newcommand{\id}{{\mathit id}}

\newcommand{\A}{{\mathcal A}}

\newcommand{\CB}{{\mathcal B}}

\def    \F      {{\mathbb F}}

\def    \R      {{\mathbb R}}

\def    \Z      {{\mathbb Z}}
\def    \N      {{\mathbb N}}

\def    \T      {{\mathbb T}}
\def    \CP     {{\mathbb C}{\mathbb P}}

\def    \12     {{\frac{1}{2}}}

\def    \HF     {\operatorname{HF}}

\def    \H      {\operatorname{H}}

\def    \CF      {\operatorname{CF}}

\def    \Fix     {\operatorname{Fix}}

\def    \bx      {\bar{x}}

\newcommand    \htop  {\operatorname{h_{\scriptscriptstyle{top}}}}

\begin{document}


\setlength{\smallskipamount}{6pt}
\setlength{\medskipamount}{10pt}
\setlength{\bigskipamount}{16pt}





\title [A generalized pseudo-rotation with positive topological entropy]{A generalized pseudo-rotation with positive topological entropy}

\author[Erman \c C\. inel\. i]{Erman \c C\. inel\. i}

\address{E\c C: Institut de Math\'ematiques de Jussieu - Paris Rive
  Gauche (IMJ-PRG), 4 place Jussieu, Boite Courrier 247, 75252 Paris
  Cedex 5, France} \email{erman.cineli@imj-prg.fr}



\date{\today}

\begin{abstract}
In this note we give examples of Hamiltonian diffeomorphisms which are on one hand dynamically complicated, for instance with positive topological entropy, and on the other hand minimal from the perspective of Floer theory. The minimality is in the sense that the barcode of the Floer complex of all iterates of these maps consists of only infinite bars. In particular, the maps have zero barcode entropy. 
\end{abstract}

\maketitle

\tableofcontents

\section{Introduction}
\label{sec:intro}

The main goal of this note is to give examples of Hamiltonian diffeomorphisms whose dynamical complexity is invisible at the level of (filtered) Floer theory. Namely, we construct  Hamiltonian diffeomorphisms which are minimal from the  perspective of Floer theory and at the same time have positive topological entropy. Here the minimality is in the sense that the barcode of the Floer complex of all iterates of these maps consists of only infinite bars. 

More precisely, let $\varphi $ be a Hamiltonian diffeomorphism of a closed symplectically aspherical symplectic manifold $(M,\omega)$. We show that the product $(M \times \T^{2n>2}, \omega \oplus \omega_{irr})$, where $(\T^{2n>2}, \omega_{irr})$ is an irrational torus (Zehnder's torus \cite{Ze}), admits a Hamiltonian diffeomorphism $\psi$ which contains on an invariant set the dynamics of $\varphi$, in particular $\htop(\psi) \geq \htop(\varphi)$, and which satisfies the following properties: The periodic points of $\psi$ are the fixed points, the fixed point set is a union of smooth submanifolds, all iterates $\psi^k$ are Morse-Bott non-degenerate and their Floer barcode consists of only infinite bars. Following \cite{AS} we call such maps generalized pseudo-rotations.

Our motivation for searching examples as above comes from the main results of \cite{CGG:Entropy, CGG:Growth} where barcode entropy was introduced. The barcode entropy is  a Floer-theoretic invariant of a Hamiltonian diffeomorhism which measures the exponential growth under iterations of the number of not-too-short bars in the barcode of the Floer complex. It was shown in \cite{CGG:Entropy} that in general the topological entropy bounds from above the barcode entropy and the two quantities are equal for Hamiltonian diffeomorphisms of closed surfaces (these results are extended to geodesic flows in \cite{GGM} and the former one to the Reeb setting in \cite{FLS}). Here we show that the latter result, the equality in dimension two, does not hold unconditionally in higher dimensions. 

In dimension two, there is a strong relationship between topological entropy and periodic point growth rate for smooth maps \cite{Ka}. Namely, by Katok's result \cite[Cor. 4.4]{Ka}, on closed surfaces the topological entropy bounds from below the exponential growth rate of periodic points. This is not necessarily true in higher dimensions. It is easy to construct diffeomorphisms, even symplectomorphisms, with no periodic points and positive topological entropy. For instance, one can product any given map with an irrational rotation of a torus to kill the periodic points while preserving the entropy. 

There cannot be such an example in the Hamiltonian case. By the Arnold's conjecture, every Hamiltonian  diffeomorphism of a closed symplectic manifold has a fixed point. Or, by the Conley conjecture \cite{Gi:CC, GG:Rev, SZ}, for a broad class of closed symplectic manifolds, which includes symplectically aspherical ones, every Hamiltonian diffeomorphism has infinitely many periodic orbits. Then, depending on the underlying symplectic manifold, one can ask for an example with finitely many periodic points (or a pseudo-rotation), or with sub-exponential periodic point growth rate, and positive topological entropy. To the best of our knowledge, there is no known example of either type. 

It is worth mentioning that pseudo-rotations can have very interesting dynamics. For instance, the complex projective space $\CP^n$ admits one which is uniquely ergodic on the complement of the fixed point set \cite{AK, LS}. However, it is not expected that a pseudo-rotation of $\CP^n$ can have positive topological entropy; see \cite[Cor.\ 2.16]{CGG:Growth} for some indirect evidence for this conjecture. We also note that a certain class of pseudo-rotations of $\CP^n$ are known to be $C^0$-rigid \cite{Br, GG:PR, JS} and hence necessarily have zero topological entropy \cite{A-Z}.

Let us turn back to the examples of this note. Our construction heavily relies on dynamical properties of the irrational flow on Zehnder's torus which forces the underlying manifold to be at least six dimensional. Another drawback of our examples is that although all periodic points are the fixed points, there are infinitely many of them (which has to be the case by the Conley conjecture).  It would be interesting to see an example of a strongly non-degenerate Hamiltonian diffeomorphism of a closed symplectic manifold where the periodic orbit growth rate is strictly less than the topological entropy.

\medskip\noindent{\bf Acknowledgements.} The author is grateful to Pierre Berger, Viktor Ginzburg, Ba\c sak G\" urel and Sobhan Seyfaddini for useful discussions. This work is supported by the ERC Starting Grant 851701.

\section{Main results}
\label{sec:cons}

In this section we explain the details of our construction and state the main properties of the resulting maps. We start by listing the ingredients:

\begin{itemize}
\item[--] a closed symplectically aspherical symplectic manifold $(M, \omega)$,
\item[--]  an irrational torus $(\T^{2n>2}, \omega_{irr})$ (see Remark \ref{rmk:torus}),
\item[--] a Morse function $\alpha \colon \R/\Z \to \R$ with exactly two critical points $\theta_1, \theta_2 \in \R/\Z$,
\item[--] a smooth function $\beta \colon \R \to \R$ which satisfies the following properties:
\subitem{--} $\beta(\alpha(\theta_i))=0$ for $i \in \{1,2\}$,
\subitem{--} there exists $\eta \in \R/\Z$ such that $\beta(\alpha(\eta))=1$. 
\end{itemize}
Out of this data we define a map
\[
\FF \colon C^{\infty}(\R/\Z \times M) \to C^{\infty}(\R/\Z \times M \times \T^{2n}) 
\] 
from one-periodic in time smooth Hamiltonians on $(M, \omega)$ to those on the product $(M \times \T^{2n}, \omega \oplus \omega_{irr})$ by setting
\[
\FF(H)(t, p, z, \theta) =  C_H \alpha(\theta) + \beta(\alpha(\theta))H(t,p)
\]
where $ p \in M$, $\theta \in \R/\Z$ is the distinguished coordinate, see Remark \ref{rmk:torus}, on the irrational torus 
\[
(z, \theta) \in \T^{2n-1} \times \R / \Z =\T^{2n}
\] 
and $C_H \in \R$ is the constant given by
\[
C_H= 2 \max \vert H  \beta' \vert.
\]

\begin{Example}
One can, for instance, set $\alpha(\theta) = \beta(\theta)= \sin (2\pi \theta)$. The critical points of $\alpha$ are $\theta_1=1/4$, $\theta_2=3/4$, where at both $\beta(\alpha(\theta))=\sin (2\pi \sin (2\pi \theta))$ vanishes. As for $\eta$ one can take any solution to $\sin(2\pi \eta) =1/4$. 
\end{Example}

\begin{Remark}
\label{rmk:iteration}
Observe that $\FF$ is $C^\infty$-continuous and it commutes with taking iterations. Namely, let us denote by $H^{\sharp k}$ the $k$-th iteration of a one periodic in-time Hamiltonian $H$, that is, $H^{\sharp k}(t,p)=kH(kt, p)$. It is easy to see that $C_{H^{\sharp k}} = k C_H$ and hence $\FF(H^{\sharp k})=\FF(H)^{\sharp k}$.
\end{Remark}

\begin{Remark}
We note that the map $\FF$ does not descend to $\Ham(M,\omega)$. In general, $\varphi_{\FF(H_1)} \neq \varphi_{\FF(H_2)}$ for Hamiltonians $H_1$, $H_2$ generating the same Hamiltonian diffeomorphism $\varphi_{H_1} = \varphi_{H_2}$.
\end{Remark}

 Next we list some properties of $\FF$ which are relevant to us. 

\begin{Theorem}
\label{thm1}
For all $H \in C^{\infty}(\R/\Z \times M)$ the image $\FF(H)$ satisfies the following properties:
\begin{enumerate}
\item $\Fix (\varphi_{\FF(H)}) = \{ M\times \T^{2n-1} \times \theta_1, M\times \T^{2n-1} \times \theta_2\}$,
\item $\varphi_{\FF(H)}$ is Morse--Bott non-degenerate,
\item $\Spec(\FF(H)) = \{  C_H \alpha(\theta_1),  C_H \alpha(\theta_2) \}$,
\item The barcode of the Floer complex  of $\FF(H)$ has no finite bars, it consists of $\dim H_*(\T^{2n-1} \times M)$ many infinite bars starting at each spectral value $C_H\alpha(\theta_i)$,
\item For all $z_0 \in \T^{2n-1}$, we have 
\[
\pi_M \circ \varphi_{\FF(H)}\circ \iota_{(z_0,\eta)} =\varphi_H
\]
where $\pi_M \colon M \times \T^{2n} \to M$ is the projection to the first factor and $\iota_{(z_0,\eta)} \colon M  \to M \times \T^{2n}$ is the embedding given by $\iota_{(z_0,\eta)}(p)= (p, z_0, \eta)$.  
\end{enumerate}
\end{Theorem}

Let us give a brief summary of how we prove Theorem \ref{thm1}.   We note that it is essentially straightforward to check that the image $\FF(H)$ satisfies these properties. What we think is new here is the observation that one can find a ``non-trivial'' Hamiltonian diffeomorphism, for instance with positive topological entropy, which has similar properties for every iteration. 

The properties \textit{(1)} and \textit{(5)} can be seen from a direct computation; and then \textit{(3)}  follows from \textit{(1)} essentially by the definition of the action spectrum. For the property \textit{(2)}, recall that a Hamiltonian diffeomorphism $\varphi$ is called Morse--Bott non-degenerate, if the fixed point set $Y$ is a union of smooth submanifolds and $\ker (D_y\varphi -\id) = T_yY$ for every $y \in Y$. In Section \ref{sec:proof}  we show that the linearization of $\varphi_{\FF(H)}$ is different than identity on the fixed point set $\Fix (\varphi_{\FF(H)})$. This is enough to conclude that $\varphi_{\FF(H)}$ is Morse--Bott non-degenerate in our case since $\Fix (\varphi_{\FF(H)})$ is a union of co-dimension one submanifolds. 

One can deduce the remaining property \textit{(4)} from \textit{(1)}, \textit{(2)} and Po\' zniak's result \cite[Cor. 3.5.4]{Po}. In Section \ref{sec:proof} we give an alternative argument using continuity of barcodes. Namely, we first observe that continuity of barcodes together with property \textit{(3)} imply that it suffices to check \textit{(4)}  at a single point in the image. As this special point we choose the image of a sufficiently small autonomous Hamiltonian and then use the fact that Morse/Floer barcodes agree for $C^2$-small autonomous Hamiltonians. 

Finally we note that the properties listed in Theorem \ref{thm1} hold for iterates $\FF(H)^{\sharp k}$ as well; see Remark \ref{rmk:iteration}. In particular, the property \textit{(5)} takes the following form: For all $H \in C^{\infty}(\R/\Z \times M)$ and $k \in \N$, we have 
\[
\pi_M \circ \varphi_{\FF(H)}^k\circ \iota_{(z_0,\eta)} =\varphi_H^k
\]
where $\pi_M$, $\iota_{(z_0,\eta)}$ are as in Theorem \ref{thm1}. As a consequence, we have $\htop(\varphi_{\FF(H)}) \geq \htop(\varphi_H)$. 

\begin{Corollary}
\label{cor}
For all $H \in C^{\infty}(\R/\Z \times M)$ and $k \in \N$, we have:
\begin{enumerate}[label=(\roman*)]
\item The periodic points of $\varphi_{\FF(H)}$ are the fixed points,
\item The fixed point set $Y$ is the union of disjoint smooth submanifolds and every iterate $\varphi_{\FF(H)}^k$ is Morse--Bott non-degenerate,
\item There exists a Morse-Bott non-degenerate and perfect function $f \colon M \times \T^{2n} \to \R$ such that the critical set of $f$ is exactly $Y$,
\item The Floer barcode of $\varphi^k_{\FF(H)}$ consists of only infinite bars and hence the barcode entropy of $\varphi_{\FF(H)}$ is zero,
\item $\htop(\varphi_{\FF(H)}) \geq \htop(\varphi_H)$. 
\end{enumerate}
\end{Corollary}

For the property \textit{(iii)} one can take, for instance, $f(p,z, \theta) =\alpha(\theta)$. The remaining properties directly follow from Theorem \ref{thm1} and Remark \ref{rmk:iteration}. 

\begin{Remark}[Zehnder's torus \cite{Ze}]
\label{rmk:torus}
Here we briefly recall the construction of an irrational symplectic form $\omega_{irr}$ on $\T^{2n>2}$. Let $(x_1, y_1, \dots ,x_n, y_n) \in (\R/\Z)^{2n} = \T^{2n>2}$ be the coordinates on $\T^{2n}$. Choose a vector $(a_1, b_1, \dots, a_{n-1}, b_{n-1}) \in \R^{2n-2}$ such that $\{a_i, b_i, 1\}$ are rationally independent and set
\[
\omega_{irr} = \sum_{i=1}^n dy_i \wedge dx_i    + \sum_{i=1}^{n-1} a_i   dx_n \wedge dy_i + \sum_{i=1}^{n-1} b_i dx_i \wedge dx_n. 
\]
Observe that the vector field $X$ that solves $\omega_{irr}(X, \cdot) =- dy_n$  is given by 
\[
X=(a_1, b_1, \dots , a_{n-1}, b_{n-1} ,1,0).
\]
It follows that every autonomous Hamiltonian on $(\T^{2n}, \omega_{irr})$ which depends only on the $y_n$-coordinate generates an irrational flow on its regular levels. This is the crucial property of the symplectic form $\omega_{irr}$. In this note we refer to the coordinate $y_n$ as the distinguished coordinate and denote it by $\theta$.  As a side remark, we note that Zehnder's torus plays an important role in Herman's  counterexample to the $C^\infty$-closing lemma \cite{He1, He2}. We refer the interested reader to \cite{HZ} for a detailed discussion of Zehnder's torus and Herman's counterexample. 
\end{Remark}

\section{Preliminaries}
\label{sec:prelim}

\subsection{Conventions and notation}
\label{sec:conventions}

Let $(M,\omega)$ be a closed symplectic manifold. We say that $(M,\omega)$ is symplectically aspherical if both $\omega(A)=0$ and $\langle c_1(TM),A \rangle=0$ for all $A \in  \pi_2(M)$. 

Let $H \colon \R/\Z \colon M \to \R$ be a one-periodic in time smooth Hamiltonian on $(M,\omega)$. The Hamiltonian vector field $X_H$ of $H$ is defined by $\omega(X_H, \cdot) = -dH$. We denote by $\varphi_H$ the time-one map of the time-dependent flow of $X_H$. Such maps form the group of Hamiltonian diffeomorphisms $\Ham(M,\omega)$ of $(M,\omega)$. For $k \in \N$, the $k$-th iterate of a Hamiltonian $H$ is defined by $H^{\sharp k}(t, p) = kH(kt, p)$. The Hamiltonian diffeomorphism $\varphi_{H^{\sharp k}}$ generated by $H^{\sharp k}$ is the $k$-th iterate of $\varphi_H$. 

Let $x \colon \R /\Z \to M$ be a contractible loop. A capping of $x$ is a map $A \colon D^2 \to M$ such that $A\vert_{\partial D^2}=x$.
The action of a Hamiltonian $H$ at a capped loop $\bx = (x, A)$ is given by
\[
\A_H(\bx) = -\int_A \omega + \int_0^1 H(t,x(t)) dt.
\]
The critical points of $A_H$ on the space of capped contractible loops are exactly the capped one-periodic orbits of $X_H$. The set of critical values $\Spec(H)$ of $\A_H$ is called the action spectrum of $H$. 

When the underlying manifold $(M,\omega)$ is symplectically aspherical, the action does not dependent on the chosen capping. Observe that in that case the action $\A_H(x)$ of a constant loop $x \colon \R /\Z \to M$  is given by $\A_H(x) =\int_0^1 H(t, x(0)) dt$.

A one-periodic orbit $x$ of $X_H$ is called non-degenerate if the linearized return map $D\varphi_H \colon T_{x(0)} M \to T_{x(0)} M$ has no eigenvalues equal to one. A Hamiltonian $H$ is called non-degenerate, if all one-periodic orbits of $X_H$ are non-degenerate. The Conley--Zehnder index $\mu(\bx)$ of a non-degenerate capped orbit $\bx$ is defined as in \cite{Sa, SZ}. In this note, $\mu$ is normalized so that $\mu(\bx)=n$ when $\bx$ is a maximum, with constant capping, of a $C^2$-small autonomous Hamiltonian. When $(M,\omega)$ is symplectic aspherical, as the action, the index does not depend on the capping either.

\subsection{Filtered Floer homology and barcodes}
\label{sec:barcodes}

In this section we recall the basics of (filtered) Floer
homology and barcodes in the symplectically aspherical setting and over $\F_2$-coefficients. Our main goal is to set our notation.  We refer the
reader to \cite{FHS, HS, MS, Sa, SZ} for a detailed
treatment of Floer homology and to \cite{PS, PRSZ, UZ} for persistence homology and barcodes. The barcodes are introduced to symplectic geometry in \cite{PS} and became a very useful tool after. The work \cite{UZ} extends the setting to general symplectic manifolds. In this note we use the version in \cite{PS, PRSZ} since we don't need the general treatment of \cite{UZ}. 

Let $H \colon \R/\Z \times M \to \R$ be a non-degenerate Hamiltonian on $(M,\omega)$ and $J$ be an $\omega$-compatible almost complex structure on the tangent space $TM$. As a vector space, the Floer chain-complex $\CF(H)$ of $H$ is given by
\[
\CF(H, J) = \bigoplus \F_2 \langle x \rangle
\]
where the sum runs through all contractible one-periodic orbits $x$ of $H$. The complex $\CF(H,J)$ is graded by the Conley-Zehnder index and the Floer differential $\partial \colon \CF_*(H,J) \to \CF_{*-1}(H, J)$ is defined as follows.  Let $x, y$ be two contractible one-periodic orbits of $H$ with $\mu(x) - \mu(y) =1$. Consider solutions $u \colon \R/\Z \times \R \to M$ to the Floer equation
\[
\partial_s u(s,t) + J (u(s,t))(\partial_t u(s,t) - X_H(u(s,t)))=0 
\]
with asymptotics $\lim_{s \to - \infty} u(s,t) =x(t)$ and $\lim_{s\to + \infty } u(s,t) =y(t)$. For a generic (time-dependent) almost complex structure $J$, the solution space, modulo the $\R$-action given by the $s$-shift, is a finite set. The coefficient of $y$ in $\partial (x)$ is given by the mod two count of this set. Moreover, for generic $J$, we have $\partial^2 =0$. 

The Floer differential $\partial$ strictly decreases the action $\A_H$ and as a result one can consider the filtered complexes $\CF^{(a,b)}(H,J)$ generated by all contractible one-periodic orbits $x$ with action $\A_H(x) \in (a,b)$. The homology $\HF(H,J)$ of the complex $\CF(H,J)$ is naturally isomorphic to the singular homology $\H(M; \F_2)$ of $M$. Unlike the the total homology $\HF(H,J)$ the filtered one $\HF^{(a,b)}(H,J)$ does depend on the Hamiltonian $H$ but not on the almost complex structure $J$. 

For all $a<b$ we have the maps $\HF^{(-\infty, a)}(H) \to \HF^{(-\infty, b)}(H)$ induced from the inclusion of chain complexes. This data is an example of a persistence module and by the Normal Form Theorem \cite[Thm. 2.1.2]{PRSZ} there is an associated interval persistence module $\CB(H)$ which is called the barcode of $H$. The barcode $\CB(H)$ is a multiset formed by some finite and half-infinite intervals. Below we give an informal description of these intervals and then list some technical properties of barcodes.  All of the facts listed below are now standard and can be found, for instance, in \cite{PRSZ}. We refer the reader to \cite[Sec. 8]{PRSZ} for a detailed treatment and further references.

The total number of half-infinite intervals in the barcode $\CB(H)$ of $H$ is equal to the dimension of the total homology $\HF(H)$, so this number does not depend on the Hamiltonian $H$. These intervals start at the spectral invariants, in a sense, they represent the homologically essential contractible one-periodic orbits of $H$. The remaining orbits are responsible for the  finite bars. More precisely, each point in the action spectrum (counted with multiplicity) appears as an end point of an interval in the barcode $\CB(H)$, hence there is finite bar if and only if the complex $\CF(H)$ has more more than $\dim \HF(H)$ generators. 

So far in this section we have considered only non-degenerate Hamiltonians. One can extend the definition of $\CB(H)$ to all, not necessarily non-degenerate, Hamiltonians $H$ by continuity. Namely, the barcode $\CB(H)$ is continuous in the Hamiltonian $H$ with respect to Hofer and spectral metrics (see \cite{KS} for the latter). In this note we only need the $C^\infty$-continuity hence we omit the definition of these metrics, and also the metric on the space of barcodes. Roughly speaking, up to finite intervals of very small length, a perturbation of the Hamiltonian results in perturbation of the end points of the intervals in the barcode.

We will use the following application of the continuity in the proof of Theorem \ref{thm1}. Suppose that we have a smooth family of, not necessarily non-degenerate, Hamiltonians $H_s$ with the property that the cardinality $\vert \Spec(H_s) \vert$ is constant in $s$. Then the continuous identification of the family $\Spec(H_s)$ also identifies the barcodes $\CB(H_s)$. One can argue these statements as follows. As mentioned above, for a non-degenerate Hamiltonian $H$, the end points of the intervals in the barcode $\CB(H)$ take values in $\Spec(H)$. The same property then holds for degenerate Hamiltonians too (by continuity of the action spectrum). Now the claim follows from the $C^\infty$-continuity of barcodes. 

Finally we note that for an autonomous Hamiltonian $H$ we have another persistence module given by the maps in singular homology $\H_*(\{f<a\}) \to \H_*(\{f<b\})$ induced from inclusion of sublevels. Let us call the associated barcode the Morse barcode of $H$. A standard but non-trivial fact is that for a $C^2$-small Hamiltonian $H$ its Morse barcode agree with the Floer theoretic one $\CB(H)$.

\section{Proof of Theorem \ref{thm1}}

\label{sec:proof}

We first prove properties \textit{(1)}--\textit{(3)} and \textit{(5)}, leaving \textit{(4)} to the end. Fix a Hamiltonian $H \in C^{\infty}(\R/ \Z \times M)$. The differential $d\FF(H)$ of $\FF(H)$ is given by
\begin{align*}
d\FF(H) (t,p, z, \theta) &= C_H \alpha'(\theta) d\theta + \beta'(\alpha(\theta)) \alpha' (\theta) H_t(p) d\theta + \beta(\alpha (\theta) ) dH_t (p) \\
&=  (C_H + \beta'(\alpha(\theta)) H_t(p)) \alpha'(\theta) d\theta + \beta(\alpha (\theta) ) dH_t (p) \\
&= g_H (t, p, \theta) \alpha'(\theta) d\theta + \beta(\alpha (\theta) ) dH_t (p).
\end{align*}
where we have set $g_H=C_H + (\beta' \circ \alpha ) H$.  Let us denote by $X$ the vector field given by $\omega_{irr}(X,\cdot) =-d\theta$. Next we write the Hamiltonian vector field $X_{\FF(H)}$ of $\FF(H)$:
\[
X_{\FF(H)} (t,p, z, \theta) = g_H(t, p, \theta)  \alpha'(\theta) X +\beta(\alpha (\theta) )  X_H (t, p).
\]
Recall that $C_H=2\max \vert \beta' H \vert $ and hence $g_H >0$.

It follows that the $\T^{2n}$ component of the Hamiltonian isotopy $\varphi_{\FF(H)}^t$ is parallel to an irrational rotation on the regular levels of $\alpha$. Then the $\theta$-coordinate of any fixed point of the time-one map $\varphi_{\FF(H)}$ is necessarily a critical point of $\alpha$. On the other hand, at the critical points $\theta_1$, $\theta_2$ of $\alpha$, the composition $\beta \circ \alpha$ vanishes (see the ingredients package in Section \ref{sec:cons}) and hence the $M$ component of the flow vanishes. This completes the proof of property \textit{(1)}. 

Let us now show that $\varphi_{\FF(H)}$ is Morse-Bott non-degenerate. Denote by $Y$ the fixed point set $\Fix(\varphi_{\FF(H)})$. It is clear that 
\[
T_yY \subset \ker (D_y\varphi_{\FF(H)}-\id) 
\]
for all $y \in Y$. Note that $Y \subset M \times \T^{2n}$ has co-dimension one and hence it suffices to show that for every $y \in Y$ there exists $v \in T_y P$ such that  $D_y\varphi(v) \neq v$.

Let $\pi (p, z, \theta) =z$ be the projection to the middle $\T^{2n-1}$ component. The composition $\pi \circ \varphi_{\FF(H)}$ has the following form:
\[
\pi \circ \varphi_{\FF(H)}(p, (z, \theta)) =z+  \bar{g}_H(p, \theta) \alpha'(\theta)\pi_* X,
\]
where $\bar{g}_H(p,\theta)= \int_0^1 g_H\circ \varphi_{\FF(H)}^t dt >0$. Since $\alpha'=0$ on $Y$, we have
\[
D_{(p,z,\theta)} \pi \circ \varphi_{\FF(H)} (\partial/ \partial \theta)=  \bar{g}_H(p, \theta) \alpha''(\theta)\pi_* X
\]
for all $(p,z,\theta) \in Y$. Furthermore, since $\alpha$ is taken to be a Morse function, the right hand side does not vanish on $Y$. We conclude that $D\varphi_{\FF(H)} (\partial/ \partial \theta) \neq \partial/ \partial \theta$ on $Y$. 

Next we focus on property \textit{(3)}. Notice that the discussion above implies that all the fixed points of $\varphi_{\FF(H)}$ are actually fixed points of the entire isotopy $\varphi_{\FF(H)}^t$. Then the action spectrum is given by the time-integral of $\FF(H)$ at the fixed points. We conclude that $\Spec(\FF(H))= \{ C_H \alpha(\theta_1), C_H \alpha (\theta_2) \}$. To show that property  \textit{(5)} holds we plug in $\theta =\eta$ to the Hamiltonian vector field $X_{\FF(H)}$ and observe that, as above, the $\T^{2n}$ component of its flow is parallel to an irrational rotation on the level $\T^{2n-1} \times \eta$, whereas the $M$ component is the Hamiltonian isotopy $\varphi_H^t$ generated by $H$ since $\beta(\alpha(\eta))=1$.

It remains to prove property  \textit{(4)}. First recall that barcodes change in a continuous manner with respect to Hofer distance, hence with respect to $C^\infty$-distance (in the Hamiltonian), and secondly that the end points of the bars take values in the action spectrum. These two properties of barcodes is enough to conclude that it suffices to show that  \textit{(4)} holds at a single non-identically zero Hamiltonian. Namely, observe that $\vert \Spec(\FF(H)) \vert =2$ unless $H \equiv 0$.  Moreover, any two such Hamiltonians can be connected by a smooth homotopy not passing through $H \equiv 0$. Since the map $\FF$ is $C^\infty$-continuous too, it follows that if property  \textit{(4)} holds at $\FF(H)$ for some non-identically zero Hamiltonian $H$, then it holds for all Hamiltonians in the image of $\FF$. To finish the proof we use the fact that filtered Floer homology of a $C^2$-small autonomous Hamiltonian is isomorphic to its Morse homology (or equivalently, the fact that the two barcodes are the same). 

More precisely, let $H\equiv C $ be the constant Hamiltonian for some $C>0$. The Hamiltonian $\FF(H)= C( 2 \max \vert \beta' \vert \alpha +  \beta \circ \alpha)$ is autonomous and hence, when $C>0$ is sufficiently small, the filtered Floer homology of $\FF(H)$ is isomorphic to its filtered Morse homology. Let us assume that this is the case and show that  the filtered Morse homology of $\FF(H)$, up to scaling by $C_H= 2C \max \vert \beta'\vert$, is isomorphic to that of the function $\alpha \colon M \times \T^{2n} \to \R$. Note that this completes the proof since the Morse barcode of $\alpha$ consists of exactly $2^{n-1} \dim H_*(M)$ infinite bars starting at each critical value $\alpha(\theta_i)$. 

We argue exactly as in the previous paragraph. Consider the family of functions $\tilde{f}_t= C( 2 \max \vert \beta' \vert \alpha +  t\beta \circ \alpha)$ for $t \in [0,1]$. We have $\tilde{f}_0=C_H \alpha$ and $\tilde{f}_1=\FF(H)$. Observe that the critical values do not change along the family $\tilde{f}_t$ and hence the filtered Morse homology (or barcodes) do not change either by continuity. We conclude that $\CB(\FF(H))$ is the same as the Morse barcode of $C_H\alpha \colon M \times \T^{2n} \to \R$.


\begin{thebibliography}{CGG22}

\bibitem[AK]{AK} D.V.  Anosov, A.B. Katok, New examples in smooth
  ergodic theory. Ergodic diffeomorphisms, (in Russian), \emph{Trudy
    Moskov.\ Mat.\ Ob\v{s}\v{c}.}, \textbf{23} (1970), 3--36.
    
    \bibitem[AS]{AS}  M.S. Atallah, E. Shelukhin, Hamiltonian no-torsion, \emph{Geom.\ Topol.},
  \textbf{27} (2023), 2833--2897. 
 
   
    \bibitem[A-Z]{A-Z} A. Avila, B. Fayad, P. Le Calvez, D. Xu, Z. Zhang,
  On mixing diffeomorphisms of the disk, \emph{Invent.\ Math.},
  \textbf{220} (2020), 673--714. 
  
    \bibitem[Br]{Br} B. Bramham, Pseudo-rotations with sufficiently
Liouvillean rotation number are $C^0$-rigid, \emph{Invent.\ Math.},
\textbf{199} (2015), 561--580. 
    
    
\bibitem[\c CGG21]{CGG:Entropy} E. \c Cineli, V.L. Ginzburg,
  B.Z. G\"urel, Topological entropy of Hamiltonian diffeomorphisms: a
  persistence homology and Floer theory perspective, Preprint
  arXiv:2111.03983.
  
  
\bibitem[\c CGG22]{CGG:Growth} E. \c Cineli, V.L. Ginzburg,
  B.Z. G\"urel, On the growth of the Floer barcode, Preprint
  arXiv:2207.03613.
  
  
\bibitem[FLS]{FLS} E. Fender, S. Lee, B, Sohn, Barcode entropy for Reeb flows on contact manifolds with Liouville fillings, Preprint arXiv:2305.04770.
  
  
\bibitem[FHS]{FHS} 
A. Floer, H. Hofer, D. Salamon,
Transversality in elliptic Morse theory for the symplectic action, 
\emph{Duke Math.\ J.},
\textbf{80} (1995), 251--292. 
  
  
  
\bibitem[Gi]{Gi:CC}
V.L. Ginzburg,
The Conley conjecture, \emph{Ann.\ of Math.} (2), \textbf{172} (2010),
1127--1180.


  \bibitem[GG17]{GG:Rev} V.L. Ginzburg, B.Z. G\"urel, Conley conjecture
  revisited, \emph{Int.\ Math.\ Res.\ Not.\ IMRN}, \textbf{2017}, doi:10.1093/imrn/rnx137.
  
  \bibitem[GG18]{GG:PR} V.L. Ginzburg, B.Z. G\"urel, Hamiltonian
pseudo-rotations of projective spaces, \emph{Invent.\ Math.},
\textbf{214} (2018), 1081--1130.
  
  
  \bibitem[GGM]{GGM} V.L. Ginzburg, B.Z. G\"urel, M. Mazzucchelli,
  Barcode entropy of geodesic flows, Preprint arXiv:2212.00943; to appear in \emph{JEMS}. 
  
\bibitem[He91a]{He1} M. Herman, Exemples de flots hamiltoniens dont aucune perturbation en topologie {$C^\infty$} n'a d'orbites p\'{e}riodiques sur un ouvert de surfaces d'\'{e}nergies, \emph{C. \ R.\ Acad. \ Sci. \ Paris S\'{e}r. \ I Math. }, \textbf{312} (1991), 989--994.
  
\bibitem[He91b]{He2} M. Herman, Diff\'{e}rentiabilit\'{e} optimale et contre-exemples \`a la fermeture en topologie {$C^\infty$} des orbites r\'{e}currentes de flots hamiltoniens, \emph{C. \ R.\ Acad. \ Sci. \ Paris S\'{e}r. \ I Math. }, \textbf{313} (1991), 49--51.

  
  
\bibitem[HS]{HS}
H. Hofer, D. Salamon, 
Floer homology and Novikov rings, in \emph{The Floer Memorial Volume},
483--524, Progr.\ Math., 133, Birkh\"auser, Basel, 1995.

\bibitem[HZ]{HZ} H. Hofer, E. Zehnder, \emph{Symplectic Invariants and
    Hamiltonian Dynamics}, Birk\"auser Verlag, Basel, 1994.
   
  
  \bibitem[JS]{JS} D. Joksimovi\'c, S. Seyfaddini, A {H}\"{o}lder type inequality for the {$C^0$} distance and {A}nosov-{K}atok pseudo-rotations, Preprint arXiv:2207.11813.
  
  \bibitem[Ka]{Ka} A. Katok, Lyapunov exponents, entropy and periodic
  orbits for diffeomorphisms, \emph{Inst.\ Hautes \'Etudes Sci.\
    Publ.\ Math.}, \textbf{51} (1980), 137--173.
    

\bibitem[KS]{KS} A. Kislev, E. Shelukhin, Bounds on spectral norms and
  applications, \emph{Geom.\ Topol.}, \textbf{25} (2021), 3257--3350.
    
 
 \bibitem[LeRS]{LS} F. Le Roux, S. Seyfaddini, The Anosov-Katok method
  and pseudo-rotations in symplectic dynamics, \emph{J. Fixed Point
    Theory Appl., Claude Viterbo's 60th birthday Festschrift.}, \textbf{24} (2022), Paper No.\ 36, 39 pp. 
    
  
\bibitem[MS]{MS}
D. McDuff, D. Salamon,
\emph{J-holomorphic Curves and Symplectic Topology}, Colloquium
publications, vol.\ 52, AMS, Providence, RI, 2012.


\bibitem[PRSZ]{PRSZ} L. Polterovich, D. Rosen, K. Samvelyan, J. Zhang,
  \emph{Topological Persistence in Geometry and Analysis}, University
  Lecture Series, vol.\ 74, Amer.\ Math.\ Soc., Providence, RI, 2020.

  \bibitem[PS]{PS} L. Polterovich, E. Shelukhin, Autonomous Hamiltonian
flows, Hofer's geometry and persistence modules, \emph{Selecta
Math.\ (N.S.)}, \textbf{22} (2016), 227--296. 

\bibitem[Po]{Po} M. Po\' zniak, \emph{Floer homology, Novikov rings and clean intersections}, Thesis, University of Warwick, 1994. 

\bibitem[Sa]{Sa}
D.A. Salamon,
Lectures on Floer homology, in \emph{Symplectic Geometry and
Topology}, IAS/Park City Math.\ Ser., vol.\ 7, Amer.\ Math.\ Soc.,
Providence, RI, 1999, 143--229.

\bibitem[SZ]{SZ}
D. Salamon, E. Zehnder,
Morse theory for periodic solutions of Hamiltonian systems and the
Maslov index, \emph{Comm.\ Pure Appl.\ Math.}, \textbf{45} (1992),
1303--1360.
  
 
\bibitem[UZ]{UZ} M. Usher, J. Zhang, Persistent homology and
Floer--Novikov theory, \emph{Geom.\ Topol.}, \textbf{20} (2016),
3333--3430. 

\bibitem[Ze]{Ze} E. Zehnder, Remarks on periodic solutions on hypersurfaces,  in \emph{Periodic Solutions of Hamiltonian Systems and Related Topics} Eds.: P. Rabinowitz et al., Reidel Publishing Co. (1987), 267--279.



\end{thebibliography}
\end{document}